\documentclass[12pt]{article}

\usepackage{amsmath,amssymb,amsthm,amscd,tabularx,array}

\newcommand{\sect}[1]{\setcounter{equation}{0}\section{#1}}
\newcommand{\subsect}[1]{\subsection{#1}}

\newfont{\extra}{msbm10 scaled\magstep1}
\newcommand{\extr}[1]{\mbox{\extra #1}}

\def\C{\extr C}
\def\K{\extr K}
\def\N{\extr N}
\def\R{\extr R}
\def\k{\kappa}
\def\kt{\tilde \kappa}

\newcommand{\RL}{\triangleright\!\!\!\blacktriangleleft}
\newcommand{\LR}{\blacktriangleright\!\!\!\triangleleft}

\begin{document}

\begin{center} 
{\LARGE{\bf{Induced representations\\[0.40cm] of quantum kinematical
algebras}}}\footnote{Talk delivered by M.A. del Olmo at the XXIII International Colloquium
on Group-Theoretical  Methods in Physics, Dubna (Russia), 31.07--05.08, 2000.}
\end{center}

\bigskip\bigskip

\begin{center} 
O. Arratia$^{*}$ and M.A. del Olmo$^{**}$ 
\end{center}

\begin{center} 
$^{*}${\sl Departamento de Matem\'atica Aplicada a la
Ingenier\'{\i}a, \\ Universidad de Valladolid E-47011, Valladolid,  Spain}\\
\medskip

$^{**}${\sl Departamento de F\'{\i}sica Te\'orica, Universidad de Valladolid,
\\ E-47011, Valladolid, Spain.}\\ 
\medskip

{e-mail:oscarr@wmatem.eis.uva.es, olmo@fta.uva.es} 
\end{center}

\vskip 0.5cm

\begin{abstract}
We  construct the induced representations of the null-plane
quantum Poincar\'e and quantum kappa Galilei algebras in $(1+1)$ dimensions. 
The induction procedure makes use of the concept of module  and is based on the existence
of a pair of Hopf algebras with a nondegenerate pairing and dual bases. 
\end{abstract}


\sect{Introduction}

Quantum kinematical algebras and groups are used for the study of
$q$--deformed symmetries of the  $q$--deformed space-time, which can be 
considered   as a non-commutative homogeneous space of the quantum kinematical 
groups \cite{podles}.

It is well known in nondeformed Lie group theory  that given a Lie group $G$ and
a  closed Lie subgroup $K$ of it,  the space of functions defined in the
homogeneous space ($\simeq G/K$) carries a representation  of $G$ induced by a
representation of $K$. Hence, symmetries and homogeneous spaces are closely related
to induced representations. On the other hand, the physical interest of the
induced representations are out of doubt \cite{wigner,olmo85}.
Thus, the study of the induced representations of 
quantum kinematical groups can be useful
for determining the behaviour of physical systems endowed with deformed
symmetries. 

In this work  we present the induced representations of
the quantum $\kt$--Galilei algebra 
and the null-plane quantum Poincar\'e algebra both in $(1+1)$ dimensions. 

The induction procedure used for us has a algebraic character since  it makes use of the
theory of modules which is, from our point of view, the appropriate tool to deal with the
algebraic structures  displayed by quantum groups and algebras
\cite{olmo98,oscar99,olmo00,ao00}. A similar method has been developed by Dobrev in
\cite{dobrev1,dobrev2} and references therein. Both procedures  deal with the dual case,
closer to the classical one, constructing representations in the algebra sector. 
Also,  one can find other   papers  
extending the induction technique to the quantum case but constructing corepresentations
 of quantum groups, i.e. representations of the coalgebra sector,  from a
mathematical perspective \cite{ibort,ciccoli98} as well as  physical
\cite{giller,maslanka,bgst98}.  

\sect{Induced representations of quantum groups}

Let $H$ be a Hopf algebra and $V$ a linear vector space over a field
$\K$ ($\R$ or
$\C$). The triplet $(V, \triangleright, H)$ 
is said to be a left $H$--module if $\alpha$ is a left action of $H$ on $V$, i.e., 
 a linear map $\alpha: H \otimes V \to V$
($ \alpha : (h \otimes v)  \mapsto\alpha(h \otimes v)\equiv h
\triangleright v$)
such that
$$
 h_1\triangleright (h_2\triangleright v)= (h_1h_2) \triangleright v, \quad
 1_H \triangleright v = v,
\qquad \forall h_1,h_2 \in H,\ \forall v \in V.
$$
Right $H$--modules can be defined in a similar way.

There are two canonical modules associated to any
pair of Hopf algebras, $H, H'$  related by a nondegenerate pairing $ \langle \,
\cdot \, , \, \cdot \, \rangle$ (under these conditions $(H, H', \langle \, \cdot \, ,
\, \cdot \, \rangle)$ will be called a nondegenerate triplet):

1) The left regular module $(H, \succ, H )$ with  action  
$$
 h_1 \succ h_2 = h_1h_2,\qquad  \forall  h_1,h_2 \in H.
$$

2) The right coregular module  $(H', \prec, H )$ with  action  defined by
$$
\langle h_2, h' \prec h_1 \rangle =
\langle h_1 \succ h_2, h' \rangle,
\qquad \forall h_1,h_2 \in  H, \quad \forall h' \in H',
$$
that  using
the coproduct in $H'$ ($\Delta(h')= h'_{(1)}\otimes h'_{(2)}$) takes the form:
$h' \prec h = \langle h, h'_{(1)} \rangle h'_{(2)}$.

The induction and coinduction algorithms of algebra representations are adapted to the
Hopf algebras as follows.
 Let $(H, H', \langle \, \cdot \, , \, \cdot \, \rangle)$ be a nondegenerate triplet 
 and  $(V, \triangleright ,K)$  a left $K$--module with $K$ a subalgebra of $H$. The
carrier space, $\K^\uparrow$, of the coinduced representation is the subspace of
$H'\otimes V$ with elements $f$ such that
 \begin{equation}\label{coinducidar}
 \langle f, k h \rangle = k \triangleright  \langle f, h \rangle, \qquad
\forall k \in K, \ \forall h \in H.
\end{equation}
The pairing  used in expression (\ref{coinducidar}) is $V$--valued and is defined  by
$\langle h' \otimes v, h \rangle= \langle h', h \rangle v$, where
$h \in H,\  h' \in H',\  v \in V$.
The action $h\triangleright f$ on the coinduced module is determined by
$$
  \langle h_1 \triangleright f, h_2 \rangle =
\langle f, h_2h_1 \rangle, \qquad \forall h_2 \in H.
$$

Let $(\K, \triangleright, K)$ be a one-dimensional coinducing module.
The carrier space of the coinduced representation is the subspace
of $H' \otimes \K \simeq H'$ composed by elements $\varphi$
verifying the equivariance condition
 $\varphi \prec k = (1 \dashv k)\varphi ,\ \forall k \in K.$
The action of $H$ on $\K^\uparrow$ induced by the action of $K$ on $\K$ is given by
$$
\langle h_2 \triangleright \varphi, h_1 \rangle =
 \langle \varphi, h_1h_2\rangle, \qquad
 \forall h_
1, h_2 \in H, \quad \forall \varphi \in \K^\uparrow ,
$$
or explicitly by 
$ h \triangleright \varphi\equiv h \succ \varphi =  \langle h, \varphi_{(2)}\rangle
\varphi_{(1)}$.

It is worthy to note that to describe the induced module  the
right, $({H'}, \prec, H)$, and left, $(H', \succ, H )$,  coregular modules are both
pertinent, the former  to determine the carrier space and the last to obtain the
induced action.

Let us consider a nondegenerate triplet $(H, {H'}, \langle \, \cdot \, , \, \cdot \,
\rangle)$ with two finite sets of generators,
$\{h_1,\ldots, h_n\}$ and $\{\varphi^1, \ldots, \varphi^n\}$, such that the families
$\{ h_l=h_1^{l_1}\cdots h_n^{l_n}\}_{l \in \N^n}$ and
$\{ \varphi^m= (\varphi^1)^{m_1} \cdots (\varphi^n)^{m_n}\}_{m \in \N^n}$ ($l=(l_1,
\ldots, l_n)$, $m=(m_1, \ldots, m_n)$) are bases of
$H$ and $H'$, respectively. 
The action on the coregular module
$({H'}, \succ, H)$ is obtained after to compute the
action of the  generators
$$
h_i \succ \varphi^j= \sum_{k \in \N^n} \alpha_{ik}^j \varphi^k,
\qquad i, j \in \{  1, 2, \ldots, n \},
$$
and to extend it to the ordered polynomial
$\varphi^j= (\varphi^1)^{j_1}\cdots (\varphi^n)^{j_n}$ using the
compatibility relation between the action
 and the algebra structure in $H'$
\begin{equation}\label{primera}
 h \succ (\varphi \psi)= (h_{(1)} \succ \varphi)
(h_{(2)} \succ \psi),\qquad
 h \succ 1_{H'} = \epsilon(h) 1_{H'}.
 \end{equation}
In order  to write  explicitly the expression of the action on
a general ordered polynomial we take into account that:

1) There is a natural
representation $\rho$, associated to  $(H,\prec,H)$, of  $H$   
$$
 [\rho(h_2)](h_1)= h_1 \prec h_2.
$$

2) The action on $({H'},\succ,H)$ can be expressed in terms of 
$\rho$  using the adjoint with respect to
$\langle \, \cdot \, , \, \cdot \, \rangle$ ($f^\dagger :H'\to H'$ is
the adjoint of $f :H\to H$ if $\langle \, h ,f^\dagger(h') \rangle= \langle \, f(h)
,h'\rangle$) defined by
\begin{equation}\label{leftcoregularaction}
 h \succ \varphi= [\rho(h)]^\dagger(\varphi).
\end{equation}

If the bases $\{ h_l\}_{l \in \N^n} $ and $\{ \varphi^m\}_{m \in \N^n}$  are dual, i.e.
$\langle h_l, \varphi^m \rangle = l! \; \delta_l^m, \  \forall l, m \in \N^n$ 
(where $l! = \prod_{i=1}^n l_i!, \ \delta_l^m= \prod_{i=1}^n \delta_{l_i}^{m_i}$),
we  define  ``multiplication'' operators $\overline{h}_i$,
$\overline{\varphi}^j$ and formal derivatives
${\partial}/{\partial h_i}$,
${\partial}/{\partial \varphi^j}$ by
$$
 \begin{array}{l}
  \overline{h}_i (h_1^{l_1} \cdots h_i^{l_i} \cdots h_n^{l_n})=
  h_1^{l_1} \cdots h_i^{l_i+1} \cdots h_n^{l_n}, \\[3mm]
 \overline{\varphi}_i
\left( (\varphi^1)^{m_1} \cdots (\varphi^i)^{m_i} \cdots (\varphi^n)^{m_n}\right) 
= (\varphi^1)^{m_1} \cdots (\varphi^i)^{m_i+1} \cdots (\varphi^n)^{m_n},\\[3mm]
\displaystyle
{\frac{\partial}{\partial h_i}} (h_1^{l_1} \cdots h_i^{l_i} \cdots h_n^{l_n}) =
 l_i \; h_1^{l_1} \cdots h_i^{l_i-1} \cdots h_n^{l_n}, \\[3mm]
\displaystyle {\frac{\partial}{\partial \varphi^i}}
\left( (\varphi^1)^{m_1} \cdots (\varphi^i)^{m_i} \cdots (\varphi^n)^{m_n}\right) 
= m_i \; (\varphi^1)^{m_1} \cdots (\varphi^i)^{m_i-1} \cdots (\varphi^n)^{m_n}.
\end{array}
$$
The adjoint operators are given by 
$
\overline{h}_i^\dagger=  {\partial}/{\partial \varphi^i}$ and 
$\overline{\varphi}^{i\dagger} =  {\partial}/{\partial h_i}$.

 \sect{Null--plane quantum Poincar\'e algebra}
 \label{pns}

The null-plane quantum deformation of the $(1+1)$ Poincar\'e algebra,
$U_{z}(\mathfrak{p}(1,1))$, is a $q$--deformed Hopf algebra that in a null-plane basis,
$\{ P_+,P_-, K\}$, has the form \cite{Bal95f} 
$$
  \begin{array}{c}
[K, P_+]= \frac{-1}{z}(e^{-2 z P_+} - 1), \qquad
[K, P_-]= -2 P_-, \qquad [P_+, P_-]= 0; \\[3mm]
\Delta P_+= P_+ \otimes 1 + 1 \otimes P_+, \qquad
\Delta X = X \otimes 1 + e^{-2 z P_+} \otimes X , \quad X \in \{P_-, K\};\\[3mm]
 \epsilon(X)= 0, \quad  X \in \{P_\pm, K\}; \\[3mm]
 S(P_+)=- P_+, \qquad S(X)=- e^{2 z P_+} X, \qquad X \in \{P_-, K\}.
  \end{array}
$$
It has also the structure of bicrossproduct 
$U_z (\mathfrak{p}(1,1))= {\cal K} \RL {\cal L}$, where $\cal K$ is a commutative and
cocommutative algebra generated by $K$, and $\cal L$ is the commutative Hopf subalgebra
of $U_z (\mathfrak{p}(1,1))$ generated by $P_+$ and $P_-$.

 The dual Hopf algebra  $F_z (P(1,1))= {\cal K}^* \LR {\cal L}^*$, where 
${\cal K}^*$ is generated by $\varphi$, and ${\cal L}^*$ by  $a_+$ and $a_-$, has the 
following structure
$$
\begin{array}{c}
 [a_+, a_-]= -2 z a_-, \qquad [a_+, \varphi]= 2 z (e^{-\varphi} -1), \qquad
 [a_-, \varphi]=0;\\[3mm]
     \Delta a_\pm = a_\pm \otimes e^{\mp 2 \varphi} + 1 \otimes a_\pm \quad, \qquad
  \Delta \varphi = \varphi \otimes 1 + 1 \otimes \varphi; \\[3mm]
\epsilon(f)= 0,\quad  f \in \{a_\pm, \varphi \}; \qquad
S(a_\pm)= - a_\pm e^{\pm \varphi}, \qquad S(\varphi )=- \varphi.
\end{array}
$$

The duality between  $U_z (\mathfrak{p}(1,1))$ and $F_z (P(1,1))$
is  explicitly given by the  pairing
$$
   \langle K^m P_-^n P_+^p, \varphi^q a_-^r a_+^s \rangle =
   m! n! p! \; \delta^m_q \delta^n_r \delta^p_s.
$$
\subsect{Coregular modules}
   \label{crlpns}

As we mentioned in the previous Section we need to know the left and the right coregular
modules, $(F_z(P(1,1)), \succ, U_z(\mathfrak{p}(1,1)))$ and $(F_z(P(1,1)), \prec,
U_z(\mathfrak{p}(1,1)))$ respectively, in order to construct the induced
representations of $U_z(\mathfrak{p}(1,1)))$.

 The structure of 
$(F_z(P(1,1)), \succ, U_z(\mathfrak{p}(1,1)))$ is given by
 \begin{equation}\label{acionregp}
 \begin{array}{l}
 K \succ (\varphi^q a_-^r a_+^s) = q \varphi^{q-1} a_-^r a_+^{s}+
   2r  \varphi^q a_-^{r} a_+^s +
   \frac{1}{z}  \varphi^q a_-^{r} a_+ [(a_+-2 z)^s -a_+^s], \\[3mm]
 P_- \succ (\varphi^q a_-^r a_+^s)  = r \varphi^q a_-^{r-1} a_+^s, \qquad
 P_+ \succ (\varphi^q a_-^r a_+^s)   = s \varphi^q a_-^r a_+^{s-1}.
 \end{array}
 \end{equation}
The following equalities are basic in the demonstration of the above  result
(\ref{acionregp})
$$
    P_-^n K=  KP_-^n + 2 n P_-^n, \qquad
    P_+^n K=  KP_+^n - n \frac{1}{z}(1- e^{-2 z P_+}) P_+^{n-1},
\qquad \forall n\in N .
$$

The structure of $(F_z(P(1,1)), \prec, U_z(\mathfrak{p}(1,1)))$ is given by
 \begin{equation}\label{moduloderecha}
 \begin{split}
  & (\varphi^q a_-^r a_+^s) \prec K = q \varphi^{q-1} a_-^r a_+^{s}, \qquad
   (\varphi^q a_-^r a_+^s)  \prec P_- =  r e^{2 \varphi} \varphi^q a_-^{r-1} a_+^s,
\\[3mm]
  & (\varphi^q a_-^r a_+^s)  \prec P_+ = -\frac{1}{2z}
         \sum_{j=1}^\infty \sum_{k=0}^j \frac{1}{j} \left(\begin{array}{c}{j}\\
{k}\end{array}\right) (-1)^k
           e^{-2 \varphi} \varphi^q a_-^r (a_+ + 2k z)^s. 
 \end{split}
 \end{equation}

The proof of (\ref{moduloderecha}) starts characterizing the module
$(U_z(\mathfrak{p}(1,1)), \succ, U_z(\mathfrak{p}(1,1)))$. For that  we take into
account the following equalities
$$
   P_- K^n = (K+2)^n P_-, \qquad
  P_+ K^n = - \frac{1}{2z} \sum_{j=1}^\infty \frac{1}{j}
(K-2j)^n (1- e^{2z P_+})^j,
  \qquad
\forall n\in N ,  
$$
that allow us to obtain easily  the explicit expression of 
$(U_z(\mathfrak{p}(1,1)), \succ, U_z(\mathfrak{p}(1,1)))$
$$
\begin{array}{l}
   K \succ K^m P_-^n P_+^p= K^{m+1} P_-^n P_+^p, \qquad
   P_- \succ K^m P_-^n P_+^p= (K+2)^{m} P_-^{n+1} P_+^p, \\[3mm]
   P_+ \succ K^m P_-^n P_+^p=   -\frac{1}{2z} \sum_{j=1}^\infty
       \frac{1}{j} (K-2j)^m P_-^n (1-e^{2z P_+})^j P_+^p. 
  \end{array}$$
The corresponding endomorphisms of $U_z(\mathfrak{p}(1,1))$ are given by
$$
      \lambda(K)= \bar{K}, \qquad
      \lambda(P_-)= \bar{P}_- e^{2\frac{\partial}{\partial K}}, \qquad
      \lambda(P_+)=  -\frac{1}{2z} \sum_{j=1}^\infty
       \frac{1}{j} e^{-2 j \frac{\partial}{\partial K}}
             (1-e^{2z \bar{P}_+})^j.
 $$
The computation of the adjoints gives
 $$
\begin{array}{lll}
 \lambda(K)^\dagger &=& \frac{\partial}{\partial \varphi}, \qquad
  \lambda(P_-)^\dagger= e^{2\bar{\varphi}}
          \frac{\partial}{\partial a_-},\\[3mm] 
     \lambda(P_+)^\dagger &=&  -\frac{1}{2z} \sum_{j=1}^\infty
       \frac{1}{j} (1-e^{2z \frac{\partial}{\partial a_+} })^j
        e^{-2 j \bar{\varphi}}=
        \frac{1}{2z} \ln\left[1- e^{-2  \bar{\varphi}}
        (1-e^{2z \frac{\partial}{\partial a_+} })\right].
   \end{array}  $$
Hence, the action on
$(F_z(P(1,1)), \prec, U_z(\mathfrak{p}(1,1)))$ is given by
\begin{equation} \label{acionregpd}
f \prec K=  \frac{\partial}{\partial \varphi} f,\quad
     f \prec P_-= e^{2\bar{\varphi}}
          \frac{\partial}{\partial a_-} f, \quad
    f \prec P_+= \frac{1}{2z} \ln\left[1- e^{-2  \bar{\varphi}}
             (1-e^{2z \frac{\partial}{\partial a_+} }) \right]\, f.
\end{equation} 
The explicit  action over the basis elements $\varphi^q a_-^r a_+^s$
(\ref{moduloderecha}) is obtained using  the series expansions of the above
expressions.

\subsect{Induced representations}

Let us consider the representation of  ${\cal L}$
\begin{equation} \label{cccp}
  1 \dashv (P_-^n P_+^p)  =  \alpha_-^n  \alpha_+^p,\qquad n,p \in  \mathbb{N}, \
\ \alpha_-, \alpha_+ \in \mathbb{C} .
\end{equation}
The  carrier space,  $\mathbb{C}^\uparrow$, of the representation
of $U_z (\mathfrak{p}(1,1))$, induced by the character  (\ref{cccp}),
 is constituted by the  elements of $F_z (P(1,1))$ having the form
$$
\phi(\varphi) e^{\alpha_- a_-} e^{ \alpha_+ a_+}.
$$
The induced representation can be translated to $\mathbb{C}[[\varphi]]$ where the
action of the  generators is
$$
\begin{array}{l}
  \phi (\varphi) \dashv K= \phi'(\varphi), \qquad
  \phi (\varphi) \dashv P_-=  \phi(\varphi)  \,  \alpha_- e^{2 \varphi},\\[3mm]
  \phi (\varphi) \dashv P_+=
   \phi (\varphi)\frac{1}{2z}\ln[ 1 - e^{-2\varphi} ( 1- e^{2 z\alpha_+})] .
\end{array}
$$

A sketch of the construction of the  representations induced by the character of
$\cal L$ (\ref{cccp}) is as follows \cite{oscar99}. To know the carrier
space of the induced representation is characterized by the equivariance condition,
which when is described in terms of the left regular module 
$(F_z(P(1,1)), \succ, U_z(\mathfrak{p}(1,1)))$ is  reduced to the equations
$$
\frac{\partial}{\partial a_-} f= \alpha_- f, \qquad
    \frac{\partial}{\partial a_+} f= \alpha_+ f.
$$
which are not really differential equations, except at the  limit 
$z \rightarrow 0$. However, their general solution is
$$
      f= \phi(\varphi) e^{\alpha_- a_-}  e^{\alpha_+ a_+},
$$
which is the same to that  obtained working formally with the derivatives.

The right regular action (\ref{acionregpd}) over $f$ gives the 
expression of the induced representation
$$
 \begin{array}{lll}
     \left[\phi(\varphi) e^{\alpha_- a_-}  e^{\alpha_+ a_+}\right] \prec K &=&
         \phi'(\varphi) e^{\alpha_- a_-}  e^{\alpha_+ a_+},\\[3mm]
     \left[\phi(\varphi) e^{\alpha_- a_-}  e^{\alpha_+ a_+}\right] \prec P_- &=&
       \phi(\varphi) e^{2 \varphi} \alpha_- e^{\alpha_- a_-}  e^{\alpha_+ a_+}, \\[3mm]
     \left[\phi(\varphi)e^{\alpha_- a_-}  e^{\alpha_+ a_+}\right]\prec P_+ &=&
       \phi(\varphi) \frac{1}{2z}
       \ln\left[1- e^{-2\varphi}(1-e^{2z \alpha_+})  \right]
         e^{\alpha_- a_-}  e^{\alpha_+ a_+}.
       \end{array}
$$
Note that in reality we have two kinds of representations labeled by the pairs
$(\alpha_+,0)$ and $(\alpha_+,1)$, respectively, since we can perform the rescaling $P_-
\to P_-/\alpha_-$ and $a_- \to \alpha_- a_-$.

Let us consider now the character of $\cal K$
\begin{equation} \label{ccclocp}
   K^n \vdash 1 = c^n,\qquad  n \in \mathbb{N},\ \  c \in \mathbb{C} .
\end{equation}
 We can construct a representation
of $U_z(\mathfrak{p}(1,1))$  whose carrier space, $\mathbb{C}^\uparrow$, 
is formed by the elements of $F_z(P(1,1))$
$$
e^{c\varphi} \phi(a_-, a_+).
$$
The  action on $\mathbb{C}^\uparrow$ can be carried to the subalgebra ${\cal L}^*$
of $F_z(P(1,1))$ obtaining
$$
  \begin{array}{lll}
    K  \vdash f(a_-, a_+) &=&  [c + 2 \bar{a}_- \frac{\partial}{\partial a_-} +
   \frac{1}{z} \bar{a}_+ (e^{-2 z \frac{\partial}{\partial a_+}}-1)]f(a_-, a_+) ,\\[3mm]
    P_\pm \vdash f(a_-, a_+) &=& \frac{\partial}{\partial a_\pm} f(a_-, a_+) .
  \end{array}
$$
Effectively, the representation induced by the character of $\cal K$ (\ref{ccclocp})
presents an equivariance condition described in terms of the left  regular 
module by the equation 
${\partial f}/{\partial \varphi} = c f$, 
whose  general solution  is
$$
         f= e^{c \varphi} \phi(a_-, a_+).
$$
The restriction of the right regular action (\ref{acionregp}) over these elements
gives the  representation 
$$
   \begin{array}{lll}
    K \succ [e^{c \varphi} \phi(a_-, a_+)] &=&
         e^{c \varphi}  \left[c  -2 \bar{a}_- \frac{\partial}{\partial a_-} +
              \frac{1}{z} \bar{a}_+
         (e^{-2z \frac{\partial}{\partial a_+}}-1)\right] \phi(a_-, a_+), \\[3mm]
   P_\mp \succ [e^{c \varphi} \phi(a_-, a_+)] &=&
         e^{c \varphi} \frac{\partial}{\partial a_\mp} \phi(a_-, a_+).
  \end{array}
$$
This representation is called ``local type'' representation because when the deformation
parameter goes to zero we recover the called local representations \cite{olmo85}. Note
that the coefficient $c$ vanishes after  the ``gauge transformation'' $K \to K - c$.

 \sect{Quantum kappa Galilei algebra}
 \label{eg}

The quantum algebra $U_{\kt}(\mathfrak{g}(1,1))$,  obtained by contraction
of the $\k$--Poincar\'e \cite{lukiersky}, is characterized by the following algebraic
structure \cite{adolfo95}:
$$
\begin{array}{c}
[H, K]= - P,  \qquad [P, K]= \frac{1}{2\kt} P^2, \qquad [H, P]=0; \\[3mm]
\Delta H = H \otimes 1 + 1 \otimes H, \qquad
\Delta X = X \otimes 1 + e^{-  \frac{1}{\kt} H} \otimes X,
\quad  X \in \{P, K \}; \\[3mm]
\epsilon(X)= 0, \quad X \in \{H, P, K\}; \\[3mm]
S(H)= -H,  \qquad   S(X)= - e^{\frac{1}{\kt} H} X,  \quad  X \in \{P, K\},
\end{array}
$$
where $\kt=\k c$, being $\k$ the deformation parameter of the above mentioned
$\k$--Poincar\'e algebra.

The dual algebra $F_{\kt}(G(1,1))$ is generated by $x,t$ and $v$, and its Hopf
structure  is 
$$
\begin{array}{c}
[t, x]= - \frac{1}{\kt} x, \qquad [x,v] = \frac{1}{2\kt} v^2, \qquad [t, v]= - 
\frac{1}{\kt} v;\\[3mm]
\Delta t = t \otimes 1 + 1 \otimes t, \quad
\Delta x = x \otimes 1 + 1 \otimes x - t \otimes v ,\quad 
\Delta v = v \otimes 1 + 1 \otimes v ; \\[3mm]
 \epsilon(f)= 0, \quad f\in \{v, t, x\}; \\[3mm]
 S(v)= -v,  \qquad  S(x)= -x - t v,  \qquad   S(t)= -t.
\end{array}
$$
The pairing between both Hopf algebras is given by
$$
\langle K^m P^n H^p, v^q x^r t^s \rangle =
m! n! p! \,  \delta^m_q \delta^n_r \delta^p_s.
$$

The action of $U_{\kt}(\mathfrak{g}(1,1))$ on
the left coregular module  
$(F_{\kt}(G(1,1)), \succ, U_{\frac{1}{2\kt}}(\mathfrak{g}(1,1)))$ is 
$$
 K \succ f= \left[\frac{\partial}{\partial v}+
\frac{1}{2\kt} \bar{x} \frac{\partial ^2}{\partial x^2} -
\bar{t} \frac{\partial}{\partial x} \right]f, \qquad
 P \succ f= \frac{\partial}{\partial x} f, \qquad
H \succ f= \frac{\partial}{\partial t} f,
$$
where $f$ is an arbitrary element of $F_{\kt}(G(1,1))$. 

The action on  the right coregular module 
$(F_{\kt}(G(1,1)), \prec, U_{\kt}(\mathfrak{g}(1,1)))$ is given by
$$
f \prec K = \frac{\partial}{\partial v} f, \quad
f \prec P = \frac{\frac{\partial}{\partial x}}{
1-\frac{\bar{v}}{2\kt}  \frac{\partial}{\partial x}} f,  \quad
f \prec H = \left[ \frac{\partial}{\partial t}  -
2\kt \ln(1-\frac{\bar{v}}{2\kt}\frac{\partial}{\partial x})\right]f.
$$

Now we can  obtain a family of representations of
$U_{\kt}(\mathfrak{g}(1,1))$ coinduced  by  the character
$$
1 \dashv P^n H^p = a^n b^p, \qquad n,p \in \N , \ \  a, b \in \mathbb{C},
$$
of the abelian subalgebra of $U_{\kt}(\mathfrak{g}(1,1))$ generated by $H$ 
and $P$, whose carrier space $\C^\uparrow$ is the set of elements of $F_{\kt}(G(1,1))$
of  the form \cite{olmo00,olmo98}
$$
\phi(v) e^{ax} e^{bt}.
$$
The action on $\C^\uparrow$ can be translated to the space of formal power series
$$
\begin{array}{lll}
\phi(v) \dashv K  & = &  \phi'(v),\qquad 
\phi(v) \dashv P  =  \phi(v) \frac{a}{1- \frac{1}{2\kt} av}, \\[3mm]
\phi(v) \dashv H  & = &  \phi(v) [b + {2\kt}\ln(1-\frac{1}{2\kt} a v)].
\end{array}
$$
The gauge transformation $H \to H-b$ allows the ``gauge equivalence" of the
representations labeled by the pair $(a,b)$ and those parameterized by  $(a,0)$.  

The ``local''
representation of
$U_{\kt}(\mathfrak{g}(1,1))$ coinduced by the character of the abelian subalgebra of
$U_{\kt}(\mathfrak{g}(1,1))$ generated by $K$
$$
K^m \vdash 1 = c^m,  \qquad m \in \N , \ \  c \in \mathbb{C},
$$
has as support the subspace  of $F_{\kt}(G(1,1))$ of
elements
$$
e^{cv} \phi(x,t).
$$
The action  of $U_{\kt}(\mathfrak{g}(1,1))$ carried to the subalgebra of formal power
series $\C[[t,x]]$ is
$$
\begin{array}{lll}
 K \vdash \phi(x,t)   & = & (c- \bar{t} \frac{\partial}{\partial x}+
\frac{1}{2\kt} \bar{x} \frac{\partial^2}{\partial x^2}) \phi(x,t),\\[3mm] 
P \vdash \phi(x,t)   & = &  \frac{\partial}{ \partial x} \phi(x,t), \qquad
 H \vdash \phi(x,t)   =   \frac{\partial}{ \partial t} \phi(x,t).
\end{array}
$$
Also here, the label $c$ can be reduced to zero.

Note that in the limit of the deformation parameter goes to zero we recover
the well know induced representations of the corresponding nondeformed Lie groups. 


{\section*{Acknowledgments}}
This work has been partially supported by DGES of the Ministerio de 
Educaci\'on  y Cultura de Espa\~na under Project PB98--0360 
and the Junta de Castilla y Le\'on (Espa\~na).


\end{document}